\documentclass[12pt]{article}
\newcommand{\doublespace}{
   \renewcommand{\baselinestretch}{1.2}
   \large\normalsize}

\def \Z{\Bbb Z}
\def \C{\Bbb C}

\def \wt{{\rm wt}}

\def \dim{{\rm dim}}
\def \Res{{\rm Res}}

\def \End{{\rm End}}
\def \Hom{{\rm Hom}}

\def \<{\langle}
\def \>{\rangle}

\def \a{\alpha }
\def \e{\epsilon }
\def \l{\lambda }
\def \L{\Lambda }

\def \pf{\noindent {\bf Proof:} \,}
\def \cg{\chi_g}
\def \cg'{\chi'_g}

\def \d{\delta}

\def \1{{\textbf 1}}
\def \h{{\frak h}}
\def \qed{{\mbox{$\square$}}}
\input amssym.def
\input amssym
\doublespace
\begin{document}
\bibliographystyle{alpha}
\newtheorem{thm}{Theorem}[section]
\newtheorem{thmn}{Theorem}
\newtheorem{prop}[thm]{Proposition}
\newtheorem{cor}[thm]{Corollary}
\newtheorem{lem}[thm]{Lemma}
\newtheorem{rem}[thm]{Remark}
\newtheorem{de}[thm]{Definition}
\newtheorem{hy}[thm]{Hypothesis}
\begin{center}
{\Large {\bf Rationality, Regularity, and $C_2$ Co-finiteness}}
\\
\vspace{0.5cm} Toshiyuki Abe\footnote{Supported by JSPS Research Fellowships 
for Young Scientists.}
\\
Department of Mathematics, Graduate School of Science, Osaka
University, Toyonaka, Osaka 560-0043, Japan
\\
{\tt sm3002at@ecs.cmc.osaka-u.ac.jp}
\\
\vspace{0.5cm} Geoffrey Buhl
\\
Department of Mathematics, University of California, Santa Cruz,
CA 95064
\\
{\tt gwbuhl@math.ucsc.edu}
\\
\vspace{0.5cm} Chongying Dong\footnote{Supported by NSF grant DMS-9987656 
and a 
research grant from the
Committee on Research, UC Santa Cruz.}
\\
Department of Mathematics, University of California, Santa Cruz,
CA 95064
\\
{\tt dong@math.ucsc.edu}
\end{center}

\begin{abstract}
We demonstrate that, for CFT vertex operator algebras, $C_2$-cofiniteness
 and rationality is equivalent to regularity.  In
addition, we show that, for $C_2$-cofinite vertex operators
algebras, irreducible weak modules are ordinary modules and 
$C_2$-cofinite, and $V_L^+$ is $C_2$-cofinite.
\end{abstract}

\section{Introduction}

One of the most important conjectures in the theory of vertex operator 
algebra is perhaps the equivalence of rationality and $C_2$-cofiniteness. 
Rationality tells us that the category of 
admissible modules is semisimple (see Section 2).  
The $C_2$ condition is slightly 
more technical and deals with the co-dimension of a certain subspace of 
$V$.  
In the case of finite dimensional Lie algebras, we have an internal 
characterization for semisimplicity.  Namely, the maximum solvable radical 
is 
zero.  This condition implies the all the modules of such a Lie algebra 
will 
be completely reducible.   For vertex operator algebras,
$C_2$-cofiniteness is the conjectured internal condition 
that implies complete reducibility of modules.  

The evidence for this conjecture is overwhelming. 
Independently, rationality 
and $C_2$-cofiniteness both imply that
the number of irreducible admissible modules for a vertex
operator algebra is finite \cite{L}, \cite{DLM2}.  Also
independently, these two notions imply that irreducible admissible
modules are irreducible ordinary modules \cite{L},\cite{DLM2}.
Rationality implies that Zhu's algebra is finite dimensional and
semisimple \cite{Z}, \cite{DLM2}, while $C_2$-cofiniteness implies that
Zhu's algebra is finite dimensional. Also well-known rational vertex
operator algebras are $C_2$-cofinite. It is not surprising that 
a lot of good results in the theory of vertex operator algebras need both 
rationality and $C_2$-cofiniteness (cf. \cite{DLM3}, \cite{DM2},\cite{Z}).

The notion of regularity is given in \cite{DLM1} to deal with weak modules
for vertex operator algebras. Regularity says that any weak module is
a direct sum of irreducible ordinary modules. It is proved in \cite{DLM1} 
that rational vertex operator algebras associated to highest 
weight modules for affine Kac-Moody algebras, Virasoro algebra, and
positive definite even lattices are regular. Based on these results
a stronger conjecture is proposed in \cite{DLM1}. That is,
rationality, $C_2$-cofiniteness, and regularity are all equivalent. 
It is proved in \cite{L} that  regularity implies $C_2$
co-finiteness.  Also, by definition regularity implies
rationality.  

In this paper we prove that $C_2$-cofiniteness and rationality together 
imply regularity. The  main idea in the proof of this result
is to find a ``lowest weight vector'' in any weak module.  In 
the case of rational affine, Virasoro, and lattice vertex
operator algebras, either the singular vectors in the Verma modules
or the lattice itself play the crucial role in the search
of lowest weight vectors \cite{DLM1}. It turns out that  the right
analogue of singular vectors in an arbitrary vertex operator algebra
is $C_2$-cofiniteness. In the case of affine or Virasoro 
vertex operator algebras, the existence of singular vectors
and $C_2$-cofiniteness are equivalent.  The PBW type of spanning 
set in any weak module for a vertex operator algebra obtained in  
\cite{Bu} 
is the key for us to obtain a lowest weight vector.

The plan for this papers is as follows.  In Section 2 we
will fix notation and give some basic definitions.  In particular we
define various notions for modules for a vertex operator algebra.
In Section 3, we review several results concerning $C_2(V).$ 
Section 4 is the proof of our main result. In Section 5
we give some additional results about $C_2$-cofinite
vertex operator algebras. In particular, we prove that
$V_L^+$ is $C_2$-cofinite for any positive definite even lattice
$L,$ extending a recent result in [Y]. 

We make the assumption that the reader is somewhat familiar with
the theory of vertex operator algebras (VOAs).  We assume the
definition of a vertex operator algebra as well as some basic
properties.

\section{Modules for vertex operator algebras}

In this section we recall the various notions of modules for a vertex
operator algebra, and we also define
the terms $C_2$-cofinite, rationality and regularity (see \cite{DLM1},
\cite{Z}).  

Throughout this paper, we will work under the assumption that $V$
is of CFT type. That is, 
$V=\bigoplus_{n\geq 0}V_n$ and $V_0=\C\1.$

\begin{de}
$V$, a VOA,  is called $C_n$-cofinite for $n\geq 2$ if $V \slash C_n(V)$ 
is
finite dimensional where $C_n(V)=\{v_{-n}w \mid v,w \in V\}.$
\end{de}

When $n=2$, the hypothesis that $V$ is $C_2$-cofinite is sometime
referred to as Zhu's finiteness condition.  This $C_2$ condition appears 
in \cite{Z}, as one of the conditions needed to prove the 
modularity of certain trace functions.  In later work \cite{GN}, it was 
shown
 that $C_2$-cofiniteness is equivalent to $C_n$-cofiniteness for all 
$n\geq 2.$ 

Rationality and regularity are two different types of complete
reducibility of VOA modules.  In order to define these terms, 
we first must describe three
different types of vertex operator algebra modules.

\begin{de}
A weak $V$ module is a vector space $M$ with a linear map 
$$
\begin{array}{ll}
Y_M: & V \rightarrow \End(M)[[z,z^{-1}]]\\
 & v \mapsto Y_M(v,z)=\sum_{n \in \Z}v_n z^{-n-1},\ \ v_n \in \End(M)
\end{array}
$$
In addition $Y_M$ satisfies the following:

1) $v_nw=0$ for $n>>0$ where $v \in V$ and $w \in M$

2) $Y_M( \1,z)=Id_M$

3) The Jacobi Identity holds:
\begin{eqnarray}
& &z_0^{-1}\d ({z_1 - z_2 \over z_0})Y_M(u,z_1)Y_M(v,z_2)-
z_0^{-1} \d ({z_2- z_1 \over -z_0})Y_M(v,z_2)Y_M(u,z_1) \nonumber \\
& &\ \ \ \ \ \ \ \ \ \ =z_2^{-1} \d ({z_1- z_0 \over 
z_2})Y_M(Y(u,z_0)v,z_2).
\end{eqnarray}
\end{de}

Of the three types of modules we mention, only weak modules have no 
grading 
assumptions.  In addition, a weak module has the minimal amount of 
structure 
needed to be a vertex algebra module.


\begin{de}
An admissible $V$ module is a weak $V$ module which carries a
$\Z_+$ grading, $M=\bigoplus_{n \in \Z_+} M(n)$, such that if $v
\in V_r$ then $v_m M(n) \subseteq M(n+r-m-1).$
\end{de}

\begin{de}
An ordinary $V$ module is a weak $V$ module which carries a $\C$
grading, $M=\bigoplus_{\l \in \C} M_{\l}$, such that:

1) $dim(M_{\l})< \infty,$

2) $M_{\l+n=0}$ for fixed $\l$ and $n<<0,$

3) $L(0)w=\l w=\wt(w) w$, for $w \in M_{\l}.$
\end{de}

An ordinary module has a grading that matches the $L(0)$ action of
the Virasoro representation as well as finite dimensional graded
pieces.  The definition of of a $\C$ grading on ordinary modules
may seem weaker than a $\Z$ grading on admissible modules.  It
turns out that the finite dimensionality of graded pieces in
ordinary modules is a strong condition, and ordinary modules are
admissible.  So we have this set of inclusions.

$$\{\mbox{ordinary modules}\} \subseteq \{\mbox{admissible modules}\}
\subseteq \{\mbox{weak modules}\}$$


\begin{de}
A vertex operator algebra is called rational if every admissible
module is a direct sum of simple admissible modules.
\end{de}

That is, a VOA is rational if there is complete reducibility of admissible 
modules. It is proved in \cite{DLM2} that 
if $V$ is rational then there are only finitely many simple admissible
modules up to isomorphism and any simple admissible module is an
ordinary module. 

\begin{de}
A vertex operator algebra is called regular if every weak module
is a direct sum of simple ordinary modules.
\end{de}

Regularity is a stronger form of complete reducibility of modules.
Again the goal of this paper is to demonstrate that under the
assumption of $C_2$-cofiniteness that rationality and regularity
are equivalent. In \cite{L}, it was shown that any regular VOA is 
$C_2$-cofinite.  So by this and the definition of rationality, regularity 
implies both rationality and $C_2$-cofiniteness.

For a vertex operator algebra $V$ and a weak $V$-module $M,$ the weight 
of the operator $v_n$ on $M$ is defined as
$\wt(v_n)=\wt(v)-n-1$ if $v$ is homogeneous. 
A vector $w\in M$ is a lowest weight vector 
if $v_n w=0$ for any homogeneous $v \in V$ and $n \in \Z$ where 
$\wt(v_n)<0$.  
The space of lowest weight vectors of a module $M$ is denoted $\Omega(M)$.  It 
is important to note that a lowest weight vector is not necessarily a 
homogeneous vector, but is has lowest weight in the sense that and operator 
of negative weight acting on a lowest weight vector will kill it.

\section{PBW type spanning set}

In this section, we will talk about the important work leading up
to the main theorem.  In particular we review a result obtained
in \cite{Bu} on a spanning set of PBW type for weak modules. 

Let $\{\bar{x}^\a\}_{\a \in I}$ be a
basis of $V \slash C_2(V)$, where $\bar{x}^\a=x^\a + C_2(V)$, and
$x^\a$ is a homogeneous vector. Set $\bar{X}=\{x^\a\}_{\a \in I}$ is a
set of elements in $V$ which are representatives of a basis for $V
\slash C_2(V)$. 
We can simplify $\bar{X}$ slightly.  The vacuum, $\1$, is not in 
$C_2(V)$.  Since the only nonzero mode of the vacuum is the identity, we 
can
discard the vacuum element from the generating set without loss.
Let us select $\bar{X}$ such that $\1 \in \bar{X}$, then let
$X= \bar{X}-\{\1\}$.  In fact, due to our assumptions that $V$ is of 
CFT type and vectors in $\bar{X}$ are homogeneous, we can set $X= \bar{X}-V_0$.

\begin{thm}
\cite{GN} $V$, a vertex operator algebra,  is spanned by
elements of the form
$$x^{1}_{-n_1} x^{2}_{-n_2} \cdots x^{k}_{-n_k} \mbox{\bf 1}$$
where $n_1>n_2> \cdots >n_k > 0$ and $x^{i} \in X$ for $1
\leq i \leq k$ .
\end{thm}

It is important to note that for this theorem, we need not assume
that $V$ is $C_2$-cofinite.  However the result is more
interesting when $V$ is $C_2$-cofinite, since the generating set
is finite.  Henceforth, we assume that $V$ is $C_2$-cofinite.

The key feature of this vertex operator algebra spanning set is that each
element satisfies a no repeat  condition.  Since in the expression of a 
spanning set element the modes are strictly decreasing, each mode
appears only once.  The next result will generalize this result so that it
applies to  modules of vertex operator algebras.  For the module spanning 
set
we will not  have a no repeat condition, but we will have a finite repeat
condition.The number of allowed repetitions will depend on how large $V
\slash C_2(V)$ is.  
 
\begin{rem}{\rm 
If $V$ is $C_2$-cofinite, then for some $N>0$,
$\bigoplus_{i\geq N} V_i \subset C_2(V)$.  In fact, $N=r+1$ where $r=\max_{x \in X} \{ \wt(x)\}$.} 
\end{rem}

Set $Q=2N-2$.  Then $Q$ will be the maximum number a times a modes
can repeat in an element of the  module spanning set. Let $W$ be a 
irreducible 
weak $V$ module, and $w \in W.$ Since $X$ is a 
finite set there is a smallest non-negative integer $L$ such that 
$x_{m}w=0$
for all $x\in X$ and $m\geq L$. 

\begin{thm}
\label{mythm} Let $V$ be a $C_2$-cofinite vertex operator
algebra, and let $W$ be an weak $V$ module generated by $w \in W$.
 Then $W$ is spanned by
elements of the form
$$x^{1}_{-n_1} x^{2}_{-n_2} \cdots x^{k}_{-n_k} w$$ where $n_1
\geq n_2 \geq \cdots \geq n_k
> -L$, and $x^i \in X$ for $1 \leq i \leq k$.  
In addition, if  $n_j>0$, then $n_j > n_{j'}$ for
$j<j'$, and if $n_j \leq 0$ then $n_j=n_{j'}$ for at most $Q-1$
indices, $j'$.
\end{thm}

These conditions say the following: For an element of the spanning
set, all the modes are decreasing and strictly less than L. If the
modes are negative then they are strictly decreasing. If the modes
are nonnegative then they are not strictly decreasing. There may
be repeats of nonnegative modes, but there are at most $Q-1$
repetitions.  Here is a sort of picture of the mode restrictions:

$$\mbox{strictly decreasing}<0 \leq \mbox{Q-1 repetitions}< L$$

The last result we need is the following.

\begin{prop}\label{prop}\cite{DLM1} Let $V$ be a rational vertex operator
algebra such that any nonzero weak $V$ module contains a simple
ordinary V submodule. Then $V$ is regular.
\end{prop}

This result is used in \cite{DLM1}  to show that rational  vertex operator 
algebras
associated the Virasoro algebra, affine Lie algebras, and lattices, 
including
the moonshine module are all regular.  In a similar fashion, we use this
Proposition to show that any $C_2$-cofinite, rational vertex operator 
algebra
is regular.

\section{Main Theorem}

In this section we first show that in any weak module
we can find a lowest weight vector.  We then use this to show the
main result; $C_2$-cofiniteness and rationality is equivalent to
regularity.  

\begin{lem}\label{l1}
\label{MT1} Let $x \in X$ and $y^1_{-m_1} y^2_{-m_2} \cdots
y^l_{-m_l} w$ be an element of the spanning set of $M$.  Then
there exist modules spanning set elements, $y^{r_1}_{-m_{r_1}}
y^{r_2}_{-m_{r_2}} \cdots y^{r_l}_{-m_{r_l}}w$, such that $$x_i
y^1_{-m_1} y^2_{-m_2} \cdots y^l_{-m_l} w=\sum_{r \in R}
c_r y^{r_1}_{-m_{r_1}} y^{r_2}_{-m_{r_2}} \cdots y^{r_l}_{-m_{r_l}}w$$
where $$\wt(x_i y^1_{-m_1} y^2_{-m_2} \cdots
y^l_{-m_l})=\wt(y^{r_1}_{-m_{r_1}} y^{r_2}_{-m_{r_2}} \cdots
y^{r_l}_{-m_{r_l}})$$ and $c_r \in \C$ for all $r \in R$, $R$ a finite 
index set.
\end{lem}

\pf 
First, we must look back to the proof of Theorem \ref{mythm} \cite{Bu}.
In the proof, three identities are used to rearrange modes to put
them in the proper form.  To show that $$\wt(x_i y^1_{-m_1}
y^2_{-m_2} \cdots y^l_{-m_l})=\wt(y^{r_1}_{-m_{r_1}}
y^{r_2}_{-m_{r_2}} \cdots y^{r_l}_{-m_{r_l}})$$ \noindent
 we will show that the identities used to rewrite the expression,
$x_i y^1_{-m_1} y^2_{-m_2} \cdots y^l_{-m_l} w$, 
as $\sum_{r \in R} c_r y^{r_1}_{-m_{r_1}}
y^{r_2}_{-m_{r_2}} \cdots y^{r_l}_{-m_{r_l}}w$ preserve the weight of
the operators.

The first of the identities, we need to look at is: 
\begin{equation}\label{id1}
[ u_{-k},
v_{-q}]=\sum\limits_{i \geq 0} {{-k} \atopwithdelims () i}
(u_{i}v)_{-k-q-i}
\end{equation}
 where $u,v \in V$ and $k,q \in \Z$. 
The second identity is
\begin{eqnarray}
(u_{-r}v)_{-q}&=&\sum\limits_{i \geq 0}(-1)^i
{{-r} \atopwithdelims () i}u_{-r-i}v_{-q+i}\nonumber\\
&&-\sum\limits_{i \geq 0}(-1)^{i-r} {{-r} \atopwithdelims ()
i}v_{-r-q-i}u_{i}\label{id}
\end{eqnarray}

\noindent where $u,v \in V$ and $r,q \in \Z$.  The third identity is more
complicated. Let $x^{1}_{-1}\cdots x^{Q}_{-1}\1= \sum_{r \in R}
x^{r_1}_{-n_{r_1}} x^{r_2}_{-n_{r_2}} \cdots
x^{r_l}_{-n_{r_l}}\1$ where $x^{i}, x^{r_t} \in X$ for $1
\leq i \leq Q$ and $1 \leq t \leq l$, $l<Q$, and $n_{r_1}
> n_{r_2} > \cdots > n_{r_l} > 0$ for fixed r.  Then
\begin{eqnarray}
\lefteqn{x_{L-k}^{Q} x_{L-k}^{{Q-1}} \cdots x_{L-k}^{1} }\nonumber \\
&=&Res_z \{ Y(\sum_{r \in R} x^{r_1}_{-n_{r_1}} x^{r_2}_{-n_{r_2}}
\cdots x^{r_l}_{-n_{r_l}}\1,z)
 z^{Q(-L-1+k)-1}\}\nonumber  \\
&&- \sum_{i=1}^{Q} \sum_{\l \in \L^i_Q} \sum_{m_{i} \geq 0}
(\prod_{j=1}^{i} x^{\l_j}_{-1-m_{\l_j}})(\prod_{j=i+1}^{Q}
x^{\bar{\l}_j}_{m_{\bar\l_j}}) \label{47s2}\\
&&- \sum_{m_{j} \geq 0, 1 \leq j \leq Q} x^{Q}_{m_{Q}}
x^{Q-1}_{m_{Q-1}}\cdots x^{1}_{m_{1}}. \label{47s3}
\end{eqnarray} 

\noindent where in (\ref{47s2}), $\sum_{j=1}^{i}(-1-m_{\l_j}) +
\sum_{j=i+1}^{Q}(m_{\bar\l_j}) = Q(L-k)$, and in (\ref{47s3}),
$\sum_{j=1}^{Q}m_{j}=Q(L-k)$ and $m_{j} \neq L-k$ for some $j$.
It is clear that the both sides of the three identities have the same 
weights.
\qed

Following \cite{L} we define $C_1(V)$ to be the subspace of $V$
spanned by $u_{-1}v, L(-1)u$ for $u,v\in \oplus_{n\geq 1}V_n.$ Since
$L(-1)u=(L(-1)u)_{-1}\1=u_{-2}\1$ we see immediately 
that $C_2(V)\subset C_1(V).$ 
$V$ is called $C_1$ cofinite if $\dim V/C_1(V)$ is finite. So if $V$ is
$C_2$-cofinite then $V$ is $C_1$-cofinite. Let $Y\subset V$ be a set
of homogeneous coset representatives of $V/C_1(V)$.

\begin{lem}\label{l2} Let $M$ be a weak module for a vertex operator 
algebra $V.$ Then we have
$$\Omega(M)=\{w\in M|y_{m}w=0, y\in Y, \wt(y_m)<0\}.$$
\end{lem}

\pf
 Let $w\in M$ such that $y_{m}w=0$ if $\wt(y_m)<0$ for $y\in Y.$ 
We prove by induction on $\wt( v)$ that $v_mw=0$ if $\wt(v_m)<0$ for
homogeneous $v\in V$.  If $v\in V_0\oplus V_1$ then $v$ is in the span
of $Y.$ The result is clear. 

Now we assume that $\wt (v)>1.$  In fact we can assume that $v\in C_1(V).$
Then $v=\sum_{i=1}^su^i_{-1}v^i+ L(-1)u$ for some
homogeneous $u^i,v^i,u\in \sum_{j\geq 1}V_j.$ Since $(L(-1)u)_m=-mu_{m+1}$ 
and $\wt (u)<\wt (v)$ it is by induction assumption that $(L(-1)u)_mw=0$ 
when
$\wt((L(-1)u)_m) < 0.$ So it is enough to show that for homogeneous
$a,b\in V,$  $(a_{-1}b)_mw=0$ if $\wt (a), \wt (b)<\wt (v)$ and 
$\wt ( (a_{-1}b)_m)<0.$ For short we set $p=\wt (a)$ and $q=\wt (b).$ Then 
$\wt (a_{-1}b)=p+q$ and $\wt ((a_{-1}b)_m)=p+q-m-1.$ So $\wt 
((a_{-1}b)_m)<0$
if and only if $m\geq p+q.$ 

Let $m\geq p+q.$ By (\ref{id}) we see that
\begin{eqnarray*}
(a_{-1}b)_{m}=\sum\limits_{i \geq 0}a_{-1-i}b_{m+i}+\sum\limits_{i \geq 0}
b_{-1+m-i}a_{i}.
\end{eqnarray*}  
Since $m\geq p+q,$ $b_{m+i}w=0$ for all $i\geq 0$ by the induction 
assumption.
Also if $i\geq p$ then $\wt( a_i)<0$ and $a_iw=0.$ So 
$$  (a_{-1}b)_{m}w=\sum_{i=0}^{p-1}b_{-1+m-i}a_{i}w.$$
By (\ref{id1}) we have
$$b_{-1+m-i}a_{i}w=a_ib_{-1+m-i}w+\sum_{t\geq 0}{-1+m-i\choose 
t}(b_ta)_{-1+m-t}w.$$
Since $i< p,$ $\wt (b_{-1+m-i})<0,$ we conclude that  $a_ib_{-1+m-i}w=0.$ 
Note that $\wt (b_ta)<p+q=\wt (v)$ for $t\geq 0$ and 
$\wt ((b_ta)_{-1+m-t})=\wt ((a_{-1}b)_m)<0.$ Again by induction 
assumption,
$(b_ta)_{-1+m-t}w=0.$ The proof is complete. 
\qed

Let $w\in M$, consider the submodule $W$ generated by $w$.  By
Theorem \ref{mythm}, $W$ is spanned by elements of the form
$$y^{1}_{-m_1} y^{2}_{-m_2} \cdots
y^{l}_{-m_l} w$$ where $m_1 \geq m_2 \geq \cdots \geq m_l
> -L$.  In addition, if  $m_j>0$, then $m_j > m_{j'}$ for
$j<j'$, and if $m_j \leq 0$ then $m_j=m_{j'}$ for at most $Q-1$
indices, $j'$.

Now these repetition restrictions allow for only finitely many
non-negative modes.  Operators with large enough positive modes 
have negative weight.  The idea here is that we can only ``push down''
$w$ so far.  So we make the following definition.


\begin{de}
\label{whatisB}
Let $B=min\{\wt(y^{1}_{-m_1} y^{2}_{-m_2} \cdots y^{l}_{-m_l})\}$
where $y^{1}_{-m_1} y^{2}_{-m_2} \cdots y^{l}_{-m_l}w$ is a non
zero spanning set element.
\end{de}

So this $B$ is the furthest we can ``push down'' $w$ without killing
$w$. Again we know that there is a minimal weight, because the
repetition restrictions only allow for finitely many positive
modes.  Next we look at how VOA spanning set elements act on module 
spanning set elements. 

\begin{lem}\label{MT2} Let $V$ be $C_2$-cofinite.  Then any weak 
$V$-module
has a nonzero lowest weight vector.
\end{lem}
\pf  
Let $v=y^1_{-m_1} \cdots y^l_{-m_l} w$ be a nonzero 
module spanning set element such that 
$$\wt(y^1_{-m_1} \cdots y^l_{-m_l})=B.$$  
We shall show that $v$ lies in $\Omega(M).$ By 
Lemma \ref{l2} we only need to prove
that $x_mv=0$ for $x\in X$ and $m\geq \wt( x).$ 

By Lemma \ref{l1} 
\begin{eqnarray*}
x_mv=\sum_{r \in R}
y^{r_1}_{-m_{r_1}} y^{r_2}_{-m_{r_2}} \cdots y^{r_l}_{-m_{r_l}}w
\end{eqnarray*}
where $\wt(y^{r_1}_{-m_{r_1}} y^{r_2}_{-m_{r_2}} \cdots
y^{r_l}_{-m_{r_l}})=\wt(x_{m-n+1} y^1_{-m_1} y^2_{-m_2} \cdots
y^l_{-m_l})<B$.  
So 
$$y^{r_1}_{-m_{r_1}}
y^{r_2}_{-m_{r_2}} \cdots y^{r_l}_{-m_{r_l}}w=0$$ for all $r \in
R$, and $x_mv=0$. 
\qed

We are now in a position to prove the main theorem of this paper. 

\begin{thm}\label{mainthm}
A vertex operator algebra, $V$, of CFT type is regular if and only
if $V$ is $C_2$-cofinite and rational.
\end{thm}

\pf  
One direction is already known, so we only need to prove that if $V$ is
rational and $C_2$ cofinite then $V$ is regular. 
From Proposition \ref{prop} it is enough to prove that any weak module
$M$ contains an irreducible ordinary module. By Lemma \ref{MT2}, 
$\Omega(M)$ is not empty. 

In order to finish the proof we need to recall the theory of associative
algebra $A(V)$ (cf. \cite{DLM2} and \cite{Z}).  
For homogeneous $u,v\in V,$ we define  products $u*v$ and $u\circ v$
 as follows:
\begin{eqnarray}\label{eqn:2.1}
& &u*v={\rm Res}_{z}\left(\frac{(1+z)^{{\rm \wt}(u)}}{z}Y(u,z)v\right)
 =\sum_{i=0}^{\infty}{{\rm \wt}(u)\choose i}u_{i-1}v\nonumber  \\
& &u\circ v={\rm Res}_{z}\left(\frac{(1+z)^{{\rm \wt}(u)}}{z^2}Y(u,z)v
\right)
 =\sum_{i=0}^{\infty}{{\rm \wt}(u)\choose i}u_{i-2}v. 
\end{eqnarray} 
Then extends (\ref{eqn:2.1}) to linear products on $V.$ Let $O(V)$ be
the linear span of $u\circ v$ for $u,v\in V.$ Set $A(V)=V/O(V).$
Then $A(V)$ is an associative algebra under
multiplication $*$ and with identity $\1+O(V)$ and central 
element
$\omega+O(V).$ Moreover, $\Omega(M)$ is an $A(V)$-module such that
$u+O(V)$ acts as $o(u)$ where $o(u)=u_{\wt (u)-1}$ if $u$ is homogeneous. 
Since our $V$ is rational, $A(V)$ is a finite dimensional
semisimple associative algebra. 

We now  pick up a simple $A(V)$-submodule $Z$ of $\Omega(M).$ Then
the $V$-submodule generated by $Z$ is an ordinary irreducible
module (see \cite{DLM2} and \cite{Z}).  
\qed

\section{$C_2$-cofinite vertex operator algebras}

In this section we study weak modules for a $C_2$-cofinite vertex
operator algebra. The result in Section 4 
on existence of lowest weight vectors in weak 
modules allow us to prove that weak modules are admissible.
We also prove that $V_L^+$ is $C_2$-cofinite for any even positive 
definite lattice $L.$

\begin{de}
Let $W$ be a weak $V$ module, then define 
$C_n(W)=\{u_{-n}w \mid u \in V, w \in W\}$.  
We say that $W$ is $C_n$-cofinite if $\dim(W \slash C_n(W)) < \infty$.
\end{de}

\begin{prop}\label{cor1}
If $V$ is $C_2$-cofinite and $W$ is an irreducible weak $V$ module, then 
$W$ is $C_2$-cofinite. 
\end{prop}

\pf
Given $w \in W$, $W$ is spanned by elements of the form 
$y^1_{-m_1} y^2_{-m_2} \cdots y^l_{-m_l}w$.  Using the mode restriction 
properties of this spanning set, there exists $N \geq 0$ such that given 
and module spanning set element $y^1_{-m_1} y^2_{-m_2} \cdots y^l_{-m_l}w$ 
with $\wt(y^1_{-m_1} y^2_{-m_2} \cdots y^l_{-m_l})>N$ then $m_1 \geq 2$.  
Thus $y^1_{-m_1} y^2_{-m_2} \cdots y^l_{-m_l}w \in C_2(W)$.
Recall Definition \ref{whatisB}:  given a module spanning set element, 
$y^1_{-m_1} y^2_{-m_2} \cdots y^l_{-m_l}w$, $B$ is the minimum weight of 
$y^1_{-m_1} y^2_{-m_2} \cdots y^l_{-m_l}$ such that the element is 
nonzero.
 Now $W \slash C_n(W)$ will be spanned by module spanning 
set elements of the form,  $y^1_{-m_1} y^2_{-m_2} \cdots y^l_{-m_l}w + 
C_2(V)$
, where $B \leq \wt(y^1_{-m_1} y^2_{-m_2} \cdots y^l_{-m_l}) \leq N$.  
This 
is a finite set because of the mode restrictions on the module spanning 
set.
\qed

Next we use the result in Proposition \ref{cor1} together with 
a recent result in \cite{Y} to prove that $V_L^+$ is $C_2$-cofinite
for any positive definite even lattice $L.$ 

We recall from [B] and [FLM] the vertex operator algebra $V_L$ associated 
to 
an even positive 
definite lattice $L.$ So $L$ is a free Abelian group of finite rank
with a positive definite $\Z$-bilinear form $(,)$ such that
$(\alpha,\alpha)\in 2\Z$ for $\alpha\in L.$ Set ${\frak 
h}=\C\otimes_{\Z}L$
and let  $\hat{\frak h}={\frak h}\otimes \C[t,t^{-1}]\oplus\C c$ be the 
corresponding
affine Lie algebra. Let $M(1)=\C[h(-n)|h\in{\frak h}, n>0]$ be the unique
irreducible module for $\hat{\frak h}$ such that $c$ acts as $1$ and 
${\frak h}\otimes t^0$ acts trivially. Then as a vector space, 
$$V_{L}=M(1)\otimes \C[L]$$
where  $\C[L]$ is the group algebra of $L$.
Then $V_L$ is a rational vertex operator algebra.

Let $\theta: V_L\to V_L$ be an order 2 automorphism such that
$$\theta(\alpha_1(-n_1)\cdots\alpha_k(-n_k)\otimes e^{\alpha})=(-1)^k
\alpha_1(-n_1)\cdots\alpha_k(-n_k)\otimes e^{-\alpha}$$
for $\alpha_i\in \h,\alpha\in L$ and $n_i<0.$ Let $V_L^{\pm}$ 
be the eigenspaces of $\theta$ with eigenvalues $\pm 1.$ Then 
$V_L^+$ is a simple vertex operator algebra and $V_L^-$ is 
an irreducible $V_L^+$-module. 

\begin{thm}\label{ths6} $V_L^+$ is $C_2$-cofinite. 
\end{thm}

\pf In the case that the rank of $L$ is one, this result has been proved 
recently in \cite{Y}. 
Let the rank of $L$ be $n.$ Then
there exists a sublattice $K$ of $L$ such that the rank $K$ is $n$ and
$K$ is a direct sum of $n$ orthogonal rank one lattice $L_1,...,L_n.$ 
Let $L=\cup_{i\in L/K}(K+\lambda_i)$ be a coset decomposition.
Then 
$$V_L=\oplus_{i\in L/K}V_{K+\lambda_i}$$ 
is a direct sum of irreducible $V_K$-modules. By
Theorems 4.4 and 6.1 of \cite{DM1}, or Theorem 6.1 of \cite{DY},
$V_{K+\lambda_i}$ is irreducible $V_K^+$-module if $2\lambda_i\not\in K$
and $V_{K+\lambda_i}$ is a direct sum of two irreducible
$V_K^+$-modules otherwise. So by Proposition \ref{cor1} it is enough
to prove that $V_K^+$ is $C_2$-cofinite.

Note that 
$$V_K=V_{L_1}\otimes \cdots \otimes 
V_{L_n}=\sum_{\e_i=\pm}V_{L_1}^{\e_1}\otimes \cdots \otimes 
V_{L_n}^{\e_n}$$
(see [FHL] for the definition of tensor product of vertex operator 
algebras) 
and 
$$V_K^+=\sum_{\e_i=\pm,\prod_{i}|\e_i|=1}V_{L_1}^{\e_1}\otimes \cdots 
\otimes 
V_{L_n}^{\e_n}$$
where $|\pm|=\pm 1.$ 
Since each $V_{L_1}^{\e_1}\otimes \cdots \otimes V_{L_n}^{\e_n}$ is
an irreducible $V_{L_1}^+\otimes\cdots\otimes V_{L_n}^+$-module, by 
Proposition
\ref{cor1} again it suffices to prove
that $V_{L_1}^+\otimes\cdots\otimes V_{L_n}^+$ is $C_2$-cofinite.

By Corollary 5.17 of \cite{Y}, each $V_{L_i}^+$ is $C_2$-cofinite. So we
are led to prove the following result: If vertex operator algebras
$V^1,...,V^n$ are $C_2$-cofinite, so is $V^1\otimes\cdots\otimes V^n.$ 

In order to see that we write $V^i=C_2(V)+W^i$ for a finite dimensional
subspace $W^i.$ It is obvious that
 $V^1\otimes\cdots V^{i-1}\otimes C_2(V^i)\otimes
V^{i+1}\otimes \cdots V^n$ is contained in $C_2(V^1\otimes \cdots \otimes
V^n).$ This shows that $V^1\otimes\
\cdots\otimes V^n= C_2(V^1\otimes \cdots \otimes
V^n)+W^1\otimes\cdots\otimes W^n.$ Thus  $V^1\otimes\cdots\otimes V^n$
is $C_2$-cofinite, as desired. \qed

\begin{cor}
Let $L=\Z \a$ such that $(\a,\a) \over 2$is prime then $V^+_L$ is regular.
\end{cor}

\pf
This follows from Theorem \ref{mainthm}, Theorem 5.17 of \cite{Y}, 
and  Theorem 4.12 
of \cite{A}, which says that $V^+_{\Z \a}$ is rational if $(\a,\a) \over 2$ 
is prime.
\qed

In order to discuss other consequences, we need a result on admissible 
modules.
\begin{lem}\label{l3} Let $V$ be $C_2$-cofinite. Then a weak module is 
admissible if and only if it is a direct sum of generalized eigenspaces 
for 
$L(0).$
\end{lem}

\pf 
First let $M=\oplus_{n=0}^{\infty}M(n)$ be an admissible $V$-module with
$M(0)\ne 0.$ Then $M(0)$ is an $A(V)$-module. Since $V$ is $C_2$-cofinite,
$A(V)$ is a finite dimensional associative algebra. 
Let $U$ be a simple 
$A(V)$-submodule $M(0).$ Then $L(0)$ acts on $U$ as a constant. As a 
result,
the $V$-module generated by $U$ is an ordinary module by Theorem 
\ref{mythm}.
Let $W$ be a maximal admissible submodule of $M$ such that $W$ is a direct
sum of generalized eigenspaces for $L(0).$ We assert that $M=W.$ 
Otherwise,
$M/W$ contains an ordinary module $\bar M/W$ for some admissible
submodule $\bar M$ of $M.$ Clearly, $\bar M$ is a sum of generalized
eigenspaces for $L(0),$ a contradiction.

Conversely, if weak module $M$ is a direct sum of generalized 
eigenspaces for $L(0).$ then it is enough to prove that 
for any $\lambda\in\C,$ the subspace 
$M^{\lambda}=\sum_{n\in\Z}M_{\lambda+n}$
is an admissible submodule where $M_{\lambda+n}$ is the
generalized eigenspace for $L(0)$ with eigenvalue $\lambda+n.$
We assume that $M^{\lambda}\ne 0.$  
Clearly, $M^{\lambda}$ is a weak submodule of $M.$ We claim that
$M_{\lambda+n}=0$ if $n$ is sufficiently small. 
For any nonzero $w\in M_{\lambda+n},$ the submodule generated by $w$ is
an admissible $V$-module by Theorem \ref{mythm}.  From the
proof of Lemma \ref{MT2}, there is a nonzero lowest weight
vector whose weight is $\lambda+m$ with  $m\leq n.$  As a result
we have a simple $A(V)$-module on which $L(0)$ acts as a scalar 
$\lambda+m.$ Since $A(V)$ has only finitely many simple modules up to 
isomorphism, we can repeat this process only a finite number times, 
proving the claim. So $M^{\lambda}=\sum_{n\geq 
N}M_{\lambda+n}$ for some $N.$ Set $M^{\lambda}(n)=M_{\lambda +N+n}$ for
$n\geq 0.$ Then $M^{\lambda}=\oplus_{n\geq 0}M^{\lambda}(n)$ is 
an admissible module. 
\qed

\begin{prop}
If $V$ is a $C_2$-cofinite vertex operator algebra, then any weak
$V$ modules is admissible. 
\end{prop}

\pf
Looking back the proof of Theorem \ref{mainthm}, all we needed to 
assume 
to show that there exists a lowest weight vector in any weak module was 
$C_2$ 
co-finiteness of $V$.  So given a weak module $M$, the space
of lowest weight vectors $\Omega(M)$ is not zero. Since $V$ is 
$C_2$-cofinite,
$A(V)$ is a finite dimensional associative algebra. Then the $A(V)$-module
$\Omega(M)$ contains a simple $A(V)$-module and the $V$-submodule 
generated by the simple $A(V)$-module is an ordinary $V$-module
by Theorem \ref{mythm}.

Let $W$ be the maximal weak submodule of $M$ which is a direct sum 
of generalized eigenspaces for $L(0).$ If $M \slash W$ is nontrivial, 
then $M \slash W$  contains a nonzero  ordinary module $ \bar M \slash W$
for some weak module $\bar M$ contained in $M.$ Clearly, $\bar M$ is
a direct sum of generalized eigenspaces for $L(0).$ This is a 
contradiction.
Thus $M=W.$ By Lemma \ref{l3}, $M$ is an admissible module.
\qed

\begin{cor}\label{cc}
If $V$ is a $C_2$-cofinite vertex operator algebra, then any
irreducible weak $V$ module is an irreducible ordinary $V$ module.
\end{cor}

\pf Since any irreducible weak module is an irreducible admissible
module, it follows from Theorem \ref{mythm} that any irreducible
admissible is ordinary. \qed

We now apply these new results to existing results about fusion rules for 
admissible modules, to obtain results about the fusion rules for weak 
modules.

\begin{de}
Let $W_i (i=1,2,3)$ be weak $V$ modules.  An intertwining operator of type 
$({W_3 \atop W_1 W_2})$ is a linear map, 
$I: W_1 \rightarrow (\Hom(W_2, W_3))[[z,z^{-1}]]$ by 
$u \mapsto I(u,z)=\sum_{\a \in \C}u_\a z^{-\a -1}$, where the following 
hold for $a \in V$, $u \in W_1$, and $v \in W_2$.

1) For all $\a$, $u_{\a+n}v=0$ for $n>>0,$

2) $I(L(-1)u,z)v=\frac{d}{dz}I(u,z)v$

3) The Jacobi Identity:
\begin{eqnarray}
z_0^{-1}\d ({z_1 - z_2 \over z_0})Y(a,z_1)I(u,z_2)-
z_0^{-1} \d ({z_2- z_1 \over -z_0})I(u,z_2)Y(a,z_1) \nonumber \\
=z_2^{-1} \d ({z_1- z_0 \over z_2})I(Y(a,z_0)u,z_2)\hspace{2cm} 
\end{eqnarray}
\end{de}

The set of intertwining operators of type $({W_3 \atop W_1 W_2})$ forms a 
vector space.  The dimension of this vector space, denoted 
$\dim({W_3 \atop W_1 W_2})$, is called fusion rule or Clebsh-Gordon 
coefficient.   We say that the fusion rules are 
finite if $\dim({W_3 \atop W_1 W_2})<\infty$ for any three irreducible 
modules $W_1, W_2, W_3$.  

\begin{cor}
Let $V$ be a $C_2$-cofinite VOA, then the fusion rules for irreducible 
weak $V$-modules are finite.
\end{cor}

\pf By Corollary \ref{cc}, any irreducible weak module is an 
ordinary module. For three irreducible weak modules $W_i (i=1,2,3)$, 
we have $W_i=\oplus_{n\geq 0}W_i(n)$ with $W_i(0)\ne 0.$ Clearly,
$\dim(W_i(0)) < \infty$.  By Proposition \ref{cor1}, we know that 
$W_i$ is $C_2$-cofinite. Then $A(W_i)$ is finite dimensional
(cf. \cite{Bu}). Here $A(W)=W/O(W)$ is an $A(V)$-bimodule
for any weak module $W$ and $O(W)$ is spanned by 
elements of the form $\Res_zY(u,z)\frac{(1+z)^{\wt (u)}}{z^2}w$
for $u\in V$ and $w\in W.$  By Proposition 2.10 of \cite{L2},
$$\dim({W_3 \atop 
W_1 W_2}) \leq \dim( \Hom_{A(V)}(A(W_1) \otimes_{A(V)} W_2(0), W_3(0))).$$
This completes the proof. 
\qed

This corollary also appears in \cite{AN} with an addition finiteness 
assumption on the module $W$.  This additional finiteness assumption can be 
removed due to results in \cite{Bu}.


\end{document}